\documentclass{amsart}
\usepackage{amsmath}
\usepackage{amssymb}
\usepackage{latexsym}

\newtheorem{thm}{Theorem}

\newtheorem{prop}{Proposition}

\def \CPb {\overline{\mathbf{CP}}{}^{2}}
\def \CP {\mathbf{CP}{}^{2}} 

\def \R {\mathbf{R}}
\def \Z {\mathbf{Z}}
\def \Sig{\Sigma}

\def \SS {S^2\times S^2}
\def \vt {\vartheta}

\def \la {\langle}
\def \ra {\rangle}

\def \b {\beta}
\def \g {\gamma}

\def \lam {\lambda}

\def \bd {\partial}
\def \x {\times}
\def \- {\!\smallsetminus\!}
\def \C {\subset}

\def \bE{\bar{E}}

\def \sign{{\text{sign}}}

\def \sw {\mathcal{SW}}

\def \DD {\Delta}

\def \SIG {\widetilde{\Sig}}
\def \CC {\widetilde{C}}
\def \RR {\widetilde{R}}

\def\spinc{spin$^{\text{c}}$}

\def\HM{H\!M_{\bullet}}
\def\wh{\widehat}

\begin{document}

\baselineskip.5cm
\title {Exotic group actions on simply connected smooth 4-manifolds}
\author[Ronald Fintushel]{Ronald Fintushel}
\address{Department of Mathematics, Michigan State University \newline
\hspace*{.375in}East Lansing, Michigan 48824}
\email{\rm{ronfint@math.msu.edu}}
\thanks{R.F. was partially supported by NSF Grant DMS-0704091, R.J.S. by NSF Grant DMS-0505080,  
and N.S. by NSF RTG Grant DMS-0353717 and by DMS-0704091.}
\author[Ronald J. Stern]{Ronald J. Stern}
\address{Department of Mathematics, University of California \newline
\hspace*{.375in}Irvine,  California 92697}
\email{\rm{rstern@uci.edu}}
\author[Nathan Sunukjian]{Nathan Sunukjian}
\address{Department of Mathematics, Michigan State University \newline
\hspace*{.375in}East Lansing, Michigan 48824}
\email{\rm{sunukjia@msu.edu}}

\dedicatory{Dedicated to Jos\'{e} Maria Montesinos on the occasion of his 65th birthday.}

\begin{abstract} We produce infinite families of exotic actions of finite cyclic groups on simply connected smooth $4$-manifolds with nontrivial Seiberg-Witten invariants.\end{abstract}
\maketitle

\section{Introduction\label{Intro}}

The goal of this paper is to exhibit infinite families of exotic actions of cyclic groups on many simply connected smooth $4$-manifolds. By {\it exotic actions} we mean smooth actions on a $4$-manifold $X$ that are equivariantly homeomorphic but not equivariantly diffeomorphic. Exotic orientation-reversing free involutions on $S^4$  were first constructed around 1980 by the first two authors \cite{invol} and can also be constructed using examples of Cappell-Shaneson \cite{CS} and later work showing that the covers of many of their manifolds are standard \cite{G1,A,G2}. In \cite{U} Ue shows that for any nontrivial finite group $G$ there is a $4$-manifold that has infinitely many free $G$-actions such that their orbit spaces are homeomorphic but mutually nondiffeomorphic. The manifolds which support Ue's exotic actions are of the form $\SS\# Z$ with $b^+(Z)>0$, and hence their Seiberg-Witten invariants vanish.

In contrast, we shall produce infinite families of finite cyclic group actions on simply connected manifolds with nontrivial Seiberg-Witten invariants.  Our theorem is:
  
\begin{thm}\label{main} Let $Y$ be a simply connected $4$-manifold with $b^+\ge1$ containing an embedded surface $\Sig$ of genus $g\ge  1$ of nonnegative self-intersection. Suppose that $\pi_1(Y\- \Sig)=\Z_d$ and that the pair $(Y,\Sig)$ has a nontrivial relative Seiberg-Witten invariant. Suppose also that $\Sig$ contains a nonseparating loop which bounds an embedded $2$-disk whose interior lies in $Y\- \Sig$. Let $X$ be the (simply connected) $d$-fold cover of $Y$ branched over $\Sig$. Then $X$ admits an infinite family of smoothly distinct but topologically equivalent actions of $\Z_d$. 
\end{thm}

As far as we know, these are the first examples of exotic orientation-preserving actions of finite cyclic groups on $4$-manifolds with nontrivial Seiberg-Witten invariants. Most of these manifolds which arise in practice are irreducible, and, in fact, if $X$ is spin with a nontrivial Seiberg-Witten invariant, $X$ must be irreducible. Construction of the group actions that we describe are obtained by  altering  branch set  data and has its origins in papers of Giffen and  
Gordon \cite{Gi,Go}.

As a simple example of our theorem, let $\Sig$ be an embedded degree $d$ curve in $\CP$.  Its complement has $\pi_1=\Z_d$ and the corresponding $d$-fold cyclic branched cover is the degree $d$ hypersurface $V_d$ in $\mathbf{CP}{}^{3}$. We can choose $\Sig$ so that it lives in a pencil and, for $d>2$, has a vanishing cycle which gives us a loop on $\Sig$ which bounds an embedded disk in its complement. That $(\CP,\Sig)$ has a nontrivial relative Seiberg-Witten invariant follows from gluing theory \cite{MST,KM}: After blowing up $d^2$ times so that the proper transform of our curve has self-intersection $0$, one can take a fiber sum with an algebraic surface containing an embedded curve of self-intersection $0$ and of the same genus as $\Sig$ to get a symplectic manifold. (See \cite{surfaces}.) Theorem~\ref{main} implies that $V_d$ admits an infinite family of topologically equivalent but smoothly distinct actions of $\Z_d$. For example, we get such a family of $\Z_4$-actions on the quartic, which is diffeomorphic to the  $K3$ surface. These examples are also discussed in the paper of H.-J. Kim \cite{K}, where, although it is not proved that the branched covers are unchanged by these operations, it is observed that the Seiberg-Witten invariants remain the same, even without the hypothesis of the theorem that there be a nonseparating loop which bounds an embedded $2$-disk.

Similarly, one can obtain an infinite family of exotic involutions on  the $K3$ surface by realizing it as the double branched cover of the sextic. For this application one needs to restate our theorem to apply to simply connected $d'$-fold branched covers where $d'$ divides $d$. This extension is, more or less, automatic, and we will not comment on it further. One can obtain an infinite family of $\Z_3$-actions on the $K3$ surface as follows. Consider a smooth embedded curve $\Sig$ in $\SS$ representing 
$3([S^2\x \{\text{pt}\}]+[ \{\text{pt}\}\x S^2])$; for example, view $\SS$ as the ruled surface $F_0$ and take $\Sig$ to be a smooth representative of the homology class of $3$ times a section plus a fiber. Then $\Sig$ has genus $4$ and $\pi_1(\SS\- N_\Sig)$ is abelian by the generalized Zariski Conjecture \cite{N}, hence 
$\pi_1(\SS\- N_\Sig)=\Z_3$. The gluing argument above implies that the relative Seiberg-Witten invariant of 
$\SS\- N_\Sig$ is nonzero, so Theorem~\ref{main} applies. In fact, using the formulas 
\[ e(X)=d\,e(Y)-(d-1)\,e(\Sig),\ \ \  \ \sign(X)= d\,\sign(Y)-\frac{(d-1)(d+1)}{3d}\Sig\cdot\Sig\]
for the euler characteristic and signature of a cyclic branched cover, one can show via a simple case-by-case analysis, 
that the only finite cyclic groups which can act on $K3$ with a smooth connected $2$-dimensional fixed point set are $\Z_2$, $\Z_3$, and $\Z_4$, and we have seen that there are infinite families of topologically equivalent but smoothly distinct examples in all these cases. 

\section{Rim surgery}

We first remind the reader of the definition of knot surgery. If $Y$ is an oriented smooth $4$-manifold containing an embedded torus $T$ of self-intersection $0$ and $K$ is a knot in $S^3$, then  {\em knot surgery} on $T$ is the result of replacing a tubular neighborhood $T\x D^2$ of $T$ with $S^1$ times the exterior $S^3\- N_K$ of the knot \cite{KL4M}:
\[ Y_K = \left( X\- (T\x D^2)\right) \cup \left(S^1\x(S^3\- N_K)\right) \]
where $\bd D^2$ is identified with a longitude $\ell_K$ of $K$. This description doesn't necessarily determine $Y_K$ up to diffeomorphism; however, when $T$ represents a non-trivial homology class in $Y$ and under reasonable hypotheses, all manifolds obtained from the same $(Y,T)$ and $K\C S^3$ will have the same Seiberg-Witten invariant:  $\sw_{Y_K}=\sw_X\cdot\DD_K(t^2)$ where where $t$ corresponds to $T$ and $\DD_K$ is the symmetrized Alexander polynomial of $K$. 

When $\Sig$ is a smoothly embedded genus $g>0$ surface in $Y$, then a relative version of knot surgery  called {\em rim surgery} \cite{surfaces} can be applied to alter the embedding type of $\Sig$. If $C$ is a homologically nontrivial loop in $\Sig$, then the preimage of $C$ under the projection of the normal circle bundle of $\Sig$ is called a {\em rim torus}. ``Rim surgery" is the result of knot surgery on a rim torus. Note that a rim torus represents a trivial homology class in $Y$ and a nontrivial homology class in $Y\setminus \Sigma$. Rim surgery replaces $C\x \bd D^2_\nu \x D^2_\delta$ with $C\x (S^3\- N_K)$ where $\bd D^2_\nu$ is the boundary circle of a normal disk to $\Sig$ and $N_K$ is a tubular neighborhood of $K$ in $S^3$. If we denote the homology class of the boundary circle of a normal disk $D^2_\delta$ to the rim torus $C\x \bd D^2_\nu$ by $\delta$,  and similarly set $\bd D^2_\nu = \nu$ and the homology class of the meridian and longitude to $K$ by $m_K$ and $\ell_K$, then the rim surgery gluing is: 
\[ \psi: C\x \bd D^2_\nu \x D^2_\delta\to S^1\x \bd(S^3\- N_K)\]
where $\psi_*([C])=[S^1]$, $\psi_*(\nu) =m_K$, and $\psi_*(\delta) =\ell_K$.
In \cite{surfaces} it is explained that this is equivalent to:
\[ (Y,\Sig_{K,C}) = (Y,\Sig) \- (C\x (I\x D^2_\nu, I\x \{ 0\}) )\cup (C\x (D^3, K')) \]
where $(D^3, K') = (S^3, K)\- (D^3, D^1)$, the knot minus a standard ball pair. This construction depends on a framing of the restriction of the normal bundle of $\Sig$ to $C\x I$ in the sense that different choices of pushoffs of $C$ to the boundary of the normal bundle may give rise to different surfaces.

In \cite{surfaces,addendum} it is shown that if both $Y$ and $Y\- \Sig$ are simply connected then $(Y,\Sig)$ and $(Y,\Sig_{K,C})$ are homeomorphic pairs, but if the self-intersection of $\Sig$ is nonnegative and if $\DD_K(t)\not\equiv 1$ then there is no self-diffeomorphism of $Y$ which throws $\Sig_{K,C}$ onto $\Sig$ provided $(Y,\Sig)$ has a nontrivial relative Seiberg-Witten invariant. In fact, under the same hypotheses, the same is true for $(Y,\Sig_{K_1,C})$ and $(Y,\Sig_{K_2,C})$ provided that the Alexander polynomials $\DD_{K_1}$ and $\DD_{K_2}$ have different sets of nonzero coefficients. Perhaps the best way to understand this (at least in the case where the self-intersection of $\Sig$ is zero) is that rim surgery multiplies the relative Seiberg-Witten invariant in the monopole Floer homology group $H\!M(\Sig\x S^1)$ for an appropriate spin$^c$-structure by $\DD_K(t^2)$. Recently, Tom Mark  \cite{M} has shown that the above result is true if the self-intersection number of $\Sig$ is greater than $2 -2g(\Sig)$. In \S4 we shall discuss this topic further.

Since a surface with a simply connected complement has no branched covers, the hypothesis that $Y\- \Sig$ is simply connected is not useful for the purpose of this paper. Kim and Ruberman \cite{K,KR} have generalized rim surgery in such a way that the condition $\pi_1(Y\- \Sig) = \Z_d$ ($d>1$) is preserved. For these purposes they used a {\em twist-rim surgery} \cite{KR} that we now describe. In rim surgery $C\x (I\x D^2_\nu, I\x \{ 0\}) \C (Y,\Sig)$ is replaced with $S^1\x (D^3, K')$. A key observation is that this last term occurs naturally in the process of spinning a knot. Given a knot $K$ in $S^3$, after removing a standard ball pair one obtains a knotted arc $K'\C D^3$. The corresponding spun 2-knot  in $S^4 = (S^1\x D^3)/\{ (t,x)\sim (t',x) \}$ (for all\ $t,t'\in S^1$ and $x\in \bd D^3$) is $S_K=(S^1\x K')/\sim$. This spun $2$-sphere $S_{K}$  naturally determines another $2$-sphere  $T_K = (S^1\x \bd D^3)/\sim$ in $S^{4}$. $T_K$ is an unknotted $2$-sphere in $S^{4}$ because it bounds the $3$-ball $\{\text{pt}\}\x D^3$. Since the  2-spheres $S_K$ and $T_K$ intersect transversely in two points, $S_K$ and $T_K$ are {\it Montesinos twins} \cite{Mo}. These twin $2$-spheres have a neighborhood $P$ in $S^{4}$ which is obtained by plumbing together two copies of $S^2\x D^2$ at two points. We call $P$ a {\it twin neighborhood}. Note that $P$ has a natural embedding in $S^4$ as the complement of the neighborhood of a standardly embedded torus in $S^4$: $S^4=P\cup (T^2\x D^2)$.

Returning to rim surgery, we identify $S^1\x (D^3, K')$ with $(S^4\- N_{T_K},  S'_K)$  where $N_{T_K}$  is a  tubular neighborhood of the $2$-sphere $T_K$, and $S'_K = S_K\cap (S^4\- N_{T_K})\cong S^1\x I$.  So rim surgery is given by the formula
\[ (Y,\Sig_{K,C}) = (Y,\Sig) \- (C\x (I\x D^2, I\x \{ 0\}) )\cup (S^4\- N_{T_K}, S'_K) \]

The process of $k$-twist-spinning a knot \cite{Z} also produces a pair of twins in $S^4$, the twist-spun knot $S_{K,k}$ and the twin $T_{K,k}$, which again arises from $\bd D^3$. (We shall give an explicit description of twist-spinning below.) The twin $T_{K,k}$ is again unknotted in $S^4$. Let $N_{T_{K,k}}$ denote a tubular neighborhood of $T_{K,k}$, then $S^4\- N_{T_{K,k}}$ is diffeomorphic to $S^1\x D^3$. One defines $k$-twist-rim surgery on $\Sig\C Y$ by
\[ (Y,\Sig_{K,C,k}) = (Y,\Sig) \- (C\x (I\x D^2, I\x \{ 0\}) )\cup (S^4\- N_{T_{K,k}}, S'_{K,k}) \]
where $S'_{K,k}=S_{K,k}\cap (S^4\- N_{T_{K,k}})$. Once again, this depends on a choice of framing for $(C\x I)\x D^2$. As we explain below, different framings may affect the value of $k$. Nonetheless, we will not further complicate matters by notationally keeping track of the framing.

The theorem of Kim and Ruberman is:

\begin{prop}[\cite{KR}] \label{KR} Let $Y$ be a simply connected smooth $4$-manifold with an embedded surface $\Sig$ of positive genus. Suppose that $\pi_1(Y\- \Sig)$ is a finite cyclic group $\Z_d$ and let $k$ be any integer relatively prime to $d$. Then for any knot $K\C S^3$ and homologically essential loop $C\C \Sig$ and for an appropriate choice of framing described below, $\pi_1(Y\- \Sig_ {K,C,k}) = \Z_d$ and, in fact, $(Y,\Sig)$ and $(Y,\Sig_ {K,C,k})$ are homemorphic as pairs. 
\end{prop}

As in the case of ordinary rim surgery, there is also a knot surgery description of $k$-twist-rim surgery. Consider the rim torus $C\x \bd D^2_\nu$ as above. Twist-rim surgery is accomplished by removing a neighborhood $C\x \bd D^2_\nu\x D^2_\delta$ of the rim torus and gluing in $S^1\x (S^3\- N_K)$ by the diffeomorphism
\[ \psi_k: C\x \bd D^2_\nu \x D^2_\delta \to S^1\x \bd(S^3\- N_K)\]
where $\psi_{k,*}([C])=k\, m_K+[S^1]$, $\psi_*[\nu] =m_K$, and $\psi_*[\delta] =\ell_K$.
The image of $\Sigma$ is now the $k$-twist rim surgered surface. Since the longitude of $K$ is identified with $\bd D^2_\delta$, we have the same (relative) Seiberg-Witten invariant as for ordinary rim surgery. This is discussed further below. 

\section{Twist-spinning and circle actions}\label{S^1}

There is a relation between twist-spinning a knot and smooth circle actions on $4$-manifolds which we shall describe in this section. Smooth circle actions on $S^4$ are completely determined by their orbit space data \cite{S^4,S^1-1,Pa}\footnote{Although these papers are set in the category of locally smooth actions, their results apply {\em{verbatim}}, with the same proofs, in the smooth category.}. The orbit space is $S^3$ or $B^3$, and in the latter case, the boundary is the image of the fixed point set and the rest of the action is free. In case the orbit space is $S^3$, the fixed point set is a pair of points, and the image of the exceptional orbits is either empty, a single arc connecting the two fixed point images, or a pair of arcs which meet only at their endpoints, the fixed point images. In case there is just one arc, its interior points all correspond to orbits with the same finite cyclic isotropy group $\Z_k$ and its endpoints to fixed points. (Thus its preimage in $S^4$ is a $2$-sphere.) We denote $S^4$ with this action by $S^4(k)$. If there are two arcs, we get a circle which contains two fixed point images splitting this circle into two arcs which correspond to finite cyclic isotropy groups of relatively prime orders. The knot type $K$ in $S^3$ which this provides is an invariant of the $S^1$-action. If the exceptional orbit types are $\Z_k$ and $\Z_d$, we denote $S^4$ with this action by $S^4(K;k,d)$. By $\bE_k$ we denote the $2$-sphere in $S^4$ consisting of the closure of the set of orbits of isotropy type $Z_k$, this is the preimage of a closed arc in $S^3$ contained in $K$. The orbit data described here completely determines smooth $S^1$-actions on $S^4$ up to equivariant diffeomorphism.

In $S^4(K;k,d)$, the $2$-spheres $\bE_k$ and $\bE_d$ form a pair of twins. The corresponding twin neighborhood is denoted $P(K;k,d)$. We will often use the notation $S^4(K;k,1)$ and $P(K;k,1)$. This gives us the $S^1$-action $S^4(k)$, but picks out a preferred set of twins in $S^4$, $\bE_k\cup \bE_1$, where $\bE_1$ is the preimage of the closed arc in $K$ labelled `$1$'. For the actions $S^4(K;k,d)$, $d\ge1$, we have 
\[ S^4(K;k,d) = P\cup (S^1\x(S^3\- N_K))\]
where $P= P(K;k,d)$ and $S^1$ acts freely in the obvious fashion on the other summand. In order to describe how these pieces are glued together, we choose bases for $H_1$ of $\bd P\cong T^3$ and $\bd(S^1\x (S^3\- N_K))$. To get such a basis for $H_1(\bd P)$ we consider the standard embedding of $P$ in $S^4$ with complement $T^2\x D^2$. (This corresponds to $K =$ unknot.)
We let $\mu_1$ be the homology class of the meridian of one of the twin two spheres, $\mu_2$ the homology class of the meridian of the other, and $\lam$ the homology class of a loop on $\bd P$ which generates $H_1(P)\cong\Z$ and which is homologically trivial in $S^4\- P$. We use the ordered basis $\{\mu_1,\mu_2,\lam\}$. For an ordered basis of $H_1(S^1\x (S^3\- N_K))$ we choose $\{ m_K, [S^1], \ell_K\}$. The gluing for $S^4(K;k,d)$, $\psi:\bd P\to \bd(S^1\x (S^3\- N_K))$ has $\psi_*$ given by the matrix
\[ A(k,d)=
\left(
\begin{array}{ccc}
k  &d  &0 \\
 -\b &\g   &0 \\
 0 &0  &1  
\end{array}
\right)\ \ \ \b\,d+\g\,k=1
\]
There is an easy description of twist-spinning in this language. If $K$ is a knot in $S^3$ then its $k$-twist-spin is $\bE_1$ in $S^4(K;k,1)$. (See e.g. \cite{Pa}.)

We describe the above surgery operations one last time in terms of this notation. If $K$ is a knot in $S^3$, let $(B^3, K')$ be $(S^3,K)$ with a trivial ball pair removed. Consider the semifree $S^1$-action on $S^4$ whose orbit space is $B^3$ with orbit map $\pi$. Then $S_K$, the spun knot obtained from $K$, is $S_K=\pi^{-1}(K')$ and its twin is $T_K=\pi^{-1}(\bd B^3)$. Rim surgery can now be described as
\[ (Y,\Sig_{K,C}) = (Y,\Sig) \- (C\x (I\x D^2, I\x \{ 0\}) )\cup (\pi^{-1}(B^3_0), \pi^{-1}(K'\cap B^3_0))\]
where $B^3_0$ is $B^3$ with an open collar of its boundary removed. Notice that from $Y$ we have removed $C\x I\x D^2 \cong S^1\x B^3$ and replaced it with $\pi^{-1}(B^3_0)\cong S^1\x B^3$, leaving the ambient space $Y$ unchanged. 

Similarly, using the $S^1$-action $S^4(K;k,1)$, the formula for $k$-twist-rim surgery becomes
\[ (Y,\Sig_{K,C,k}) = (Y,\Sig) \- (C\x (I\x D^2, I\x \{ 0\}) )\cup (S^4\- N(\bE_k), \bE_1') \]
where $N(\bE_k)$ is an $S^1$-equivariant tubular neighborhood and $\bE_1' =\bE_1\cap (S^4\- N(\bE_k))$.
Because $S^4(K;k,1)=P\cup (S^1\x (S^3\- N_K))$, we can express $k$-twist-rim surgery in terms of surgery on the torus $R=C\x \bd D^2_\nu$:
\[ (Y,\Sig_{K,C,k}) = (Y\- ( R\x D^2_\delta) ,\Sig) \cup_\phi (S^1\x  (S^3\- N_K),\emptyset) \]
where $\phi_*$ is the composition $\psi_*\circ \iota$ with $\iota_*(C)=\mu_1$, $\iota_*(\nu)=\mu_2$ and $\iota_*(\delta)=\lam$. We are identifying $\mu_1$ with a normal circle to $\bE_k$ and $\mu_2$ with a normal circle to $\bE_1$. (See also \cite{Pl} where this is described in slightly different notation.) The matrix giving our gluing is the matrix $A(k,1)$ defined above.

As we have pointed out, our construction depends on a choice of framing for the restriction of the normal bundle of $\Sig$ to $C\x I$ or equivalently of the rim torus $C\x \bd D^2_\nu$. If one pushoff $C'$ of $C$ gives rise to the gluing above, then any other framing comes from replacing $C'$ by $C'+r\nu$. Thus it corresponds to the gluing matrix $A( k+r,1)$. Thus $k$-twist-rim surgery with respect to the first framing is $(k+r)$-twist-rim surgery with respect to the second framing. This brings us to the choice of framing in the Kim-Ruberman Theorem. We need to choose a framing so that $C\x \{\text{pt}\}$ is nullhomologous in $Y\- \Sig$. Because $H_1(Y\- \Sig)=\Z_d$ is generated by $\nu$, different choices of acceptable framings differ by integer multiples of $d\,\nu$. Thus $k$-twist-rim surgery gets turned into into $(k+rd)$-twist-rim surgery, preserving the hypothesis that $k$ and $d$ should be relatively prime.

\section{Branch sets and relative Seiberg-Witten invariants}

Fix an integer $d>1$. The $\Z_d$-actions which we construct will be $d$-fold cyclic branched covers of smoothly embedded surfaces in smooth $4$-manifolds. Let $Y$ be a simply connected oriented smooth $4$-manifold with $b^+(Y)\ge 1$ containing a smoothly embedded surface $\Sig$ of genus $g\ge 1$ and self-intersection $\Sig\cdot\Sig=n \ge 0$ such that $\pi_1(Y\- \Sig)=\Z_d$. 

Choose a homologically essential loop $C$ on $\Sig$ and for $k\ge 1$ relatively prime to $d$ perform $k$-twist-rim surgery on $\Sig$ using the rim torus corresponding to $C$ and a knot $K$ in $S^3$ to obtain a surface $\Sig_{K,C,k}$. We now fix $d$ and $k$ and use the shorthand $\Sig_K =\Sig_{K,C,k}$. It follows from the result of Kim and Ruberman, Proposition~\ref{KR}, that the pairs $(Y,\Sig)$ and $(Y,\Sig_K)$ are homeomorphic.  We say that surfaces $\Sig$, $\Sig'$ in $Y$ are {\em smoothly (resp. topologically) equivalent} if there is a diffeomorphism (resp. homeomorphism) of pairs $(Y,\Sig)\cong (Y,\Sig')$.

We employ a simple trick to reduce to the situation where the self-intersection of the surface is $0$. Blow up $\Sig\cdot\Sig=n$ times to get $\wh{Y}=Y\# n\,\CPb$, and let $\wh{\Sig}$ be the blown up surface, which has self-intersection $0$. If the surfaces $\Sig_1$, $\Sig_2$ are smoothly equivalent in $Y$, then $\wh{\Sig}_1$ and $\wh{\Sig}_2$ will be smoothly equivalent in $\wh{Y}$. Furthermore, $\wh{\Sig}_{K,C,k}=\wh{\Sig_{K,C,k}}$.
Thus we may assume that $\Sig\cdot\Sig=0$. There is a complete proof in \cite{surfaces}
that when the genus of $\Sigma$ is $1$, if rim surgery is performed using knots with distinct Alexander polynomials then one obtains smoothly distinct embedded surfaces, and the same holds for twist-rim surgery. Hence we may assume that $g\ge 2$.

We need to describe the Seiberg-Witten invariant $SW_{(Y|\Sig)}$ as defined in \cite{KM,KM2}. This is the Seiberg-Witten invariant of $Y\- N(\Sig)$ obtained from \spinc-structures $\mathfrak{s}$ on $Y$  which satisfy $\la c_1(\mathfrak{s}), \Sig\ra=2g-2$. On $\Sig\x S^1$ there is a unique \spinc-structure ${\mathfrak{s}}_{g-1}$ 
which is pulled back from a \spinc-structure on $\Sig$ and satisfies $\la c_1(\mathfrak{s}_{g-1}), \Sig\ra=2g-2$. As explained in \cite{KM}, systems of local coefficients for monopole Floer homology correspond to $1$-cycles on $\Sig\x S^1$, and up to isomorphism the groups depend on their homology classes
$\eta\in H_1(\Sig\x S^1;\R)$. The related monopole Floer homology groups with local coefficients are
$\HM(\Sig\x S^1| \Sig;\Gamma_{\eta})=\R$.
In fact, as is pointed out in \cite{KM2}, if we take a product metric on $\Sig\x S^1$ where the metric on $\Sig$ has constant negative curvature, then there is a unique nondegenerate solution of the Seiberg-Witten equations on $\Sig\x S^1$. This gives rise to a distinguished generator of each of the groups $\HM(\Sig\x S^1| \Sig;\Gamma_{\eta})$ and thus they can all be identified. 

In order to get an invariant of the pair $(Y,\Sig)$ we need to consider the relative Seiberg-Witten invariant of $Y\- N(\Sig)$. Relative Seiberg-Witten invariants of a $4$-manifold with boundary take their values in various Floer homology groups of the boundary. In our case, all the groups can be identified as pointed out above to obtain the relative invariant described below.

Let $W= Y\- N(\Sig)$, and assume that $W$ inherits an orientation and homology orientation from $(Y,\Sig)$. The space $\mathcal{B}(W;[a_0])$ of pairs $(A,\Phi)$ consisting of a \spinc-connection and spinor which limit to the unique equivalence class $[a_0]$ of solutions of the Seiberg-Witten equations for the \spinc-structure ${\mathfrak{s}}_{g-1}$ on $\Sig\x S^1$ splits into path components, and each path component $z$ determines a \spinc-structure $\mathfrak{s}_{W,z}$. The moduli space of solutions to the Seiberg-Witten equations on $W^*$, i.e. $W$ with a cylindrical end, splits along these path components, 
$M(W^*;[a_0])=\coprod M_z(W^*;[a_0])$, and similarly for the whole configuration space, $\mathcal{B}(W^*;[a_0])=\coprod \mathcal{B}_z(W^*;[a_0])$. The set of all path components, $\pi_0(\mathcal{B}(W^*;[a_0]))$, is a principal homogeneous space for $H^2(W,\Sig\x S^1;\Z)$.

For each such path component $z$, any pair $(A,\Phi)$ representing $z$, and any class $\nu\in H_2(W, \Sig\x S^1;\R)$, the integral over $\nu$ of the curvature $F_{A^t}$ of the connection $A^t$, induced on the determinant line of the spinor bundle, depends only on $z$ and $\nu$. 
We thus have a relative Seiberg-Witten invariant defined by
\[SW_{(Y|\Sig)}: H_2(W, \Sig\x S^1;\R)\to \R\]
\[ SW_{(Y|\Sig)}(\nu)= \sum_z m_W(z) \exp(\frac{i}{2\pi}\int_{\nu} F_{A_z^t})\]
where the sum is taken over $z\in \pi_0(\mathcal{B}(W^*;[a_0]))$.  (It is shown in \cite{KM} that for an appropriate perturbation of the Seiberg-Witten equations, only finitely many such $z$ admit solutions of the Seiberg-Witten equations.) The coefficient $m_W(z)$ denotes the count with signs of solutions in $M_z(W^*;[a_0])$ in case this moduli space is $0$-dimensional; $m_W(z)$ is $0$ otherwise. No assumption on $b^+(W)$ is necessary for the definition of the invariant $SW_{(Y|\Sig)}$. (See \cite[\S 3.9]{KM}.)

The proof of the knot surgery theorem \cite{KL4M} tells us that $SW_{(Y|\Sig_K)}$ is obtained from 
$SW_{(Y|\Sig)}$ by multiplying by $\DD_K(t)$, the symmetrized Alexander polynomial of $K$. We now explain this.  Write $\DD_K(t)= \sum_{j=-d}^d c_jt^j$, and let $\rho\in H^2(W,\Sig\x S^1;\Z)$ be the Poincar\'e dual of the rim torus $R$ corresponding to the loop $C$ on $\Sig$.  Then, recalling that $\pi_0(\mathcal{B}(W^*;[a_0]))$ is a principal homogeneous  space for $H^2(W,\Sig\x S^1;\Z)$, for each $z$ such that $m_W(z)\ne 0$, we have $z+j\rho\in \pi_0(\mathcal{B}(W^*;[a_0]))$, $j\in\Z$. Because the calculation for twist-rim surgery is the same as for rim surgery, the knot surgery theorem gives
\begin{multline*} SW_{(Y|\Sig_{K,C,k})}(\nu) =  SW_{(Y|\Sig)}(\nu)\cdot \DD_K(\rho^2) =\\ 
      =\sum_{z,j} m_W(z) c_j \exp\big{(}2j\la\rho,\nu\ra+\frac{i}{2\pi}\int_{\nu} F_{A_z^t}\big{)}= \\
        =\sum_{z,j} m_{W_K}(z+j\rho) \exp\big{(}\frac{i}{2\pi}\int_{\nu} F_{A_{z+j\rho}^t}\big{)}
\end{multline*}
where $W_K = Y\- (\Sig_K\x D^2)$, and we are identifying $H^2(W_K,\Sig_K\x S^1;\Z)$ with the group $H^2(W,\Sig\x S^1;\Z)$
using the canonical isomorphism described in \cite{surfaces}. Note that this formula asserts that $m_{W_K}(z+j\rho)= c_j\,m_W(z)$ and that $\frac{i}{2\pi}\int_{\nu} F_{A_{z+j\rho}^t} = 2j\la\rho,\nu\ra+\frac{i}{2\pi}\int_{\nu} F_{A_z^t}$.

We would like to be able to conclude that if $\DD_{K_1}(t)\ne\DD_{K_2}(t)$ then rim (or twist-rim) surgery using these two knots results in smoothly inequivalent surfaces in $Y$. However, all that we are currently able to say is the following. (Compare \cite{addendum}.) Let ${\mathcal{S}}_{W_K}=\{z\in \pi_0(\mathcal{B}(W_K^*;[a_0]))\mid m_{W_K}(z)\ne 0\}$.

\begin{prop} \label{SS}
If $\Sig_{K_1}$ and $\Sig_{K_2}$ are smoothly equivalent, there is 
an automorphism of $H^2(W,\Sig\x S^1;\Z)$ sending ${\mathcal{S}}_{W_{K_1}}$ to ${\mathcal{S}}_{W_{K_2}}$ and preserving the coefficients $m_{W_{K_i}}(z)$. 
\end{prop}
\begin{proof} If $\Sig_K$ and $\Sig_{K'}$ are smoothly equivalent in $Y$, then $\wh{\Sig}_K$ and $\wh{\Sig}_{K'}$ are smoothly equivalent in $\wh{Y}$. The proposition now follows because relative Seiberg-Witten invariants are invariants of smooth equivalence of surfaces. (Cf. \cite{addendum}.)
\end{proof}

\section{Cyclic group actions: Equivariant rim surgery}

Let $Y$ be a simply connected smooth $4$-manifold, with an embedded surface $\Sig$ of genus $g\ge 1$   whose self-intersection number is nonegative and such that $\pi_1(Y\- \Sig) = \Z_d$. Let $C$ be a nonseparating loop on $\Sig$ which bounds a disk in $Y\- \Sig$, for example, $C$ could be a vanishing cycle.

Let $K$ be a knot in $S^3$, and for an integer $k$ relatively prime to $d$ perform $k$-twist-rim surgery on $\Sig$ using the loop $C$, and let $\Sig_K=\Sig_{K,C,k}\C Y$. Proposition~\ref{KR} implies that the surfaces $\Sig$ and $\Sig_K$ are topologically equivalent in $Y$.   Let $W$ and $W_K$ be the complements of tubular neighborhoods of $\Sig$ and $\Sig_K$ in $Y$. As in the previous section, blow up $\Sig\cdot\Sig$ times to obtain $\wh{Y}$, $\wh{\Sig}$, and  $\wh{\Sig}_K$. Note that the blowup of $\Sig_K$ is the same as the result of twist-rim surgery ${\wh{\Sig}}_K$ on  ${\wh{\Sig}}$.

\begin{prop} \label{3} Let $K$ and $K'$ be two knots in  $S^3$ and suppose that their Alexander polynomials have different (unordered) sets of nontrivial coefficients. Also suppose that $SW_{(Y|\Sig)}\ne 0$. 
Then the surfaces $\Sig_K$ and $\Sig_{K'}$ are smoothly inequivalent in $Y$. In particular, the $\Z_d$-actions on the $d$-fold cyclic covers of $Y$ branched over $\Sig_K$ and $\Sig_{K'}$ are equivariantly homeomorphic but not equivariantly diffeomorphic. \end{prop}

\begin{proof} The hypothesis implies via the knot surgery formula that ${\mathcal{S}}_{W_{K}}\ne {\mathcal{S}}_{W_{K'}}$. 
The first part of the proposition now follows from Proposition~\ref{SS}. The second part of the proposition follows from Proposition~\ref{KR} and the fact that an equivariant diffeomorphism induces a diffeomorphism of orbit spaces preserving the fixed point image.
\end{proof}

We now need to show that the cyclic branched covers in question are diffeomorphic. Let $X_K$ be the $d$-fold cyclic cover of $Y$ branched over $\Sig_K$. We have seen that
\[ (Y,\Sig_K) = (Y,\Sig) \- (C\x (I\x D^2, I\x \{ 0\}))\cup (S^4(K;k,1)\- N(\bE_k), \bE_1') \]

Alternatively, we have the rim surgery description:
\[ (Y,\Sig_K)= (Y,\Sig) \- ( R\x D^2_\delta) \cup_\phi (S^1\x  (S^3\- N_K)) \]
where $S^4(K;k,1)\- P(K;k,1)\cong S^1\x (S^3\- N_K)$, 

Let $\vt: (X,\SIG)\to (Y,\Sig)$ and $\vt_K:(X_K,\SIG_K)\to (Y,\Sig_K)$ be  the branched covers. The loop $C$ lifts to a loop $\CC$ on $\SIG$ and the rim torus $R$ similarly lifts to the rim torus $\RR$ associated to $\SIG$ and $\CC$. 	  
The manifold $X_K$ is obtained from $X$ by replacing $\RR\x D^2_\delta$ with the $d$-fold cover of 	$S^4(K;k,1)\- P(K;k,1)$. According to \cite{Pa} (see also \cite{Pl}), the branched cover of $S^4(K;k,1)$ branched over $\bar{E}_1$ is $S^4(K;k,d)$, and the branching locus is the $2$-sphere $\bar{E}_d$. The deck transformations of this branched cover are generated by the action of $e^{2\pi\,i/d}\in S^1$ contained in the circle action; so the branched covering map $S^4(K;k,d) \to S^4(K;k,1)$ sends $\bE_k$ to $\bE_k$.
Thus we can see that to obtain $X_K$, we replace $\RR\x D^2_\delta$ with $S^4(K;k,d)\- P(K;k,d)$ which is in turn diffeomorphic to $S^1\x (S^3\- N_K)$.
\[ X_K= X\- (\RR\x D^2_\delta)\cup_{\tilde{\phi}} (S^1\x(S^3\- N_K)) \]
where $\tilde{\phi}$ is given by the matrix $A(k,d)$ when bases are chosen as in \S~\ref{S^1}. And again as 
in that section, we may redescribe $X_K$ as
\[ X_K= (X\- (\CC\x I\x D^2_\nu)) \cup (S^4(K;k,d)\- N(\bE_k))\]

\begin{prop} If $C$ bounds an embedded disk in $Y\- \Sig$ then $X_K$ is diffeomorphic to $X$.
\end{prop}
\begin{proof} If $C$ bounds an embedded disk in $Y\- \Sig$, in the cover this means that $\CC$ bounds a disk in $X\- \SIG$. (In fact $\CC$ bounds $d$ such disks with disjoint interiors.) 
Hence a pushoff $\CC\x \{\text{pt}\}$ bounds an embedded disk in  $X\- (\CC\x I\x D^2_\nu)$.  The union of a regular neighborhood $U$ of this disk with $\CC\x I\x D^2_\nu$ is the result of attaching a $2$-handle to $S^1\x B^3$ along $S^1\x \{\text{pt}\}$. This is the $4$-ball, $B^4$.

The $S^1\x B^3$ in question is $\CC\x I\x D^2_\nu$ and the rim torus is $\RR = \CC\x \{\text{pt}\}\x \bd D^2_\nu$ in the boundary $S^1\x S^2 = \CC\x D^2\cup_{\RR}\CC\x D^2$. Attaching the $2$-handle corresponds to surgery on $\CC\x\{\text{pt}\}$; so $\bd B^4 = \CC\x D^2\cup_{\RR} S^1\x D^2$ where the gluing takes some pushoff of $\CC$ to $\bd D^2$. If $\CC'$ is a preferred pushoff of $\CC$, i.e. it is nullhomologous in $X\- \SIG$, then our gluing takes $[\CC'] + r[\bd D^2_\nu]$ to $[\bd D^2]=0$.

Thus the rim torus $\RR$ is a standard unknotted torus in $S^3=\bd B^4$, and, if we take the union of $B^4=U\cup (\CC\x I\x D^2_\nu)$ with another copy of $B^4$, we get $S^4=P\cup_{\RR\x\bd D^2_\delta} \RR\x D^2_\delta$. After the handle addition, the standard basis $\{\mu_1,\mu_2,\lam\}$ of $H_1(\bd P)$ is identified with $\{\CC'+r\nu, \nu, \delta\}$ in $H_1(\RR\x\bd D^2_\delta)$.

Let $V=\CC\x I\x D^2_\nu\-  (\RR\x D^2_\delta)$ which is diffeomorphic $S^1$ times the standard cobordism from a torus to a $2$-sphere obtained by attaching a $2$-handle.
In $X_K$, the $4$-ball
$ U\cup (\CC\x I\x D^2_\nu) =  U\cup V \cup (\RR\x D^2_\delta)$
 is replaced by 
 \[ U\cup V \cup (S^4(K;k,d)\- P(K;k,d)) = U\cup V \cup_{\tilde{\phi}} (S^1\x (S^3\- N_K))\] 
 
 However using obvious notation, $S^4(K;k,d) =P\cup_{A(k,d)} (S^1\x (S^3\- N_K))$; so 
 \begin{multline*} (B^4\cup U\cup V) \cup_{\tilde{\phi}} (S^1\x (S^3\- N_K)) = \\
 (S^4\- {\text{Nbd}}(T^2_{\text{std}})\cup_{\tilde{\phi}} (S^1\x (S^3\- N_K)) = P\cup_A (S^1\x (S^3\- N_K))\cong S^4 
 \end{multline*}
 because 
 \[ A= A(k,d)\circ \left(
\begin{array}{ccc}
1  &0  &0 \\
r &1   &0 \\
 0 &0  &1  
\end{array}
\right)
=
\left(
\begin{array}{ccc}
k+rd &d  &0 \\
-\b+r\g &\g   &0 \\
 0 &0  &1  
\end{array}
\right)
= A(k+rd,d)
 \]
Hence  $U\cup V \cup_{\tilde{\phi}} (S^1\x (S^3\- N_K)) \cong B^4$. It follows that $X_K$ is diffeomorphic to $X$.

\end{proof}

We may now complete the proof of Theorem~\ref{main}.

\begin{proof}[of Theorem~\ref{main}] Fix a positive integer $k$ relatively prime to $d$ and a nonseparating simple closed curve $C$ on $\Sig$ such that $C$ bounds an embedded $2$-disk $D$ in $Y\- \Sig$. Let $\{K_i\}_{i=1}^{\infty}$ be a family of knots in $S^3$ whose Alexander polynomials have pairwise different sets of nonzero coefficients. 
It then follows from Propositions~\ref{KR} and \ref{3} that the surfaces $\Sig_{K,C,k}$ obtained from $\Sig$ by $k$ twist-rim surgery are topologically equivalent but smoothly distinct. Of course this means that their corresponding branched covers give $Z_d$-actions which are equivariantly homeomorphic but not equivariantly diffeomorphic. Furthermore, because each of these branched covers $X_{K_i}$ is obtained from $X$ by removing a $4$-ball and then replacing it with another $4$-ball, each $X_{K_i}$ is, in fact, diffeomorphic to $X$.
\end{proof}

Note that the construction used in the proof can be viewed as an {\it equivariant rim surgery}.  In fact we could have presented the construction in this manner. However, it has been convenient to phrase our arguments in the language of circle actions in order to more easily identify the gluing diffeomorphisms and to more clearly see that the  construction will not change $X$ as long as $C$ bounds an embedded disk in the complement of $\Sig$.

\section{Final comments}

As we have shown above, Theorem~\ref{main} applies widely. Many smooth $4$-manifolds are constructed as branched $\Z_{d}$-covers and, with mild conditions on the branch set, they thus have infinite families of exotic actions of $\Z_{d}$. In most cases these manifolds are irreducible. All these actions are  nontrivial on homology. (Because otherwise $e(X)=e(Y)$, which implies $e(Y)=e(\Sig)$. But $X$ is simply connected, so this implies $e(X)=2$, which is ruled out if $X$ has a nontrivial Seiberg-Witten invariant.) It remains an interesting question to determine if there are simply connected $4$-manifolds with exotic actions of cyclic groups $\Z_{d}$ ($d >2$) that induce the identity on homology. This is of particular interest for the $K3$ surface. Also, it is still an interesting problem to determine exotic free group actions on a fixed smooth $4$-manifold with a nontrivial Seiberg-Witten invariant. All these questions are in the realm of seeking general rigidity or uniqueness properties in dimension $4$.

\end{document}